\documentclass[reqno,showkeys]{amsart}
\usepackage{amsfonts}
\usepackage{amssymb}
\usepackage{amsthm}
\usepackage{upref}
\makeatletter
 
\def\LaTeX{\leavevmode L\raise.42ex
    \hbox{\kern-.3em\size{\sf@size}{0pt}\selectfont A}\kern-.15em\TeX}
 \makeatother

 \newcommand{\supp}{\operatorname{supp}}
  \newcommand{\e}{\eqref}
\newcommand{\ri}{\rightarrow}
\newcommand{\q}{\quad}
\newcommand{\ii}{\infty}

   \newcommand{\sgn}{\operatorname{sgn}}

\numberwithin{equation}{section}

\newtheorem{lemma}{Lemma}[section]
\newtheorem{theorem}[lemma]{Theorem} 

\newtheorem{proposition}[lemma]{Proposition}

\theoremstyle{definition}

\newtheorem{counterexample}[lemma]{Counterexample}

\newtheorem{assumption}[lemma]{Assumption}

\theoremstyle{remark}

\newcommand{\dist}{\operatorname{dist}}

\newenvironment{||}{\|}

\newenvironment{pf}{\begin{proof}}{\end{proof}}

\def\qqq{\mathrel{\subset\mkern-15mu\lower.38ex\hbox{${\scriptscriptstyle\rightarrow}$}}}
\let\cal\mathcal

\let\Bbb\mathbb

\begin{document}  
 
\title{Exponential decay of eigenfunctions of first order systems}
\author{ D. R. Yafaev}
\address{ IRMAR, Universit\'{e} de Rennes I\\ Campus de
  Beaulieu, 35042 Rennes Cedex, FRANCE}
\email{yafaev@univ-rennes1.fr}
\thanks{Research partially supported by SPECT}

\begin{abstract}
For first order systems,
we obtain an efficient bound on the exponential decay of an eigenfunction in terms of the distance between the corresponding eigenvalue and the essential spectrum. As an example, the Dirac operator is considered.
\end{abstract}
\maketitle

\thispagestyle{empty}

%***********************************************************
\section{Introduction}
%***********************************************************
%{\bf 1.1. }
There are two different methods for proving  exponential decay of eigenfunctions $\psi$   for second order selfadjoint elliptic operators $H$. The first of them is due to Shnol' \cite{Sh} who proved that an eigenfunction corresponding to an eigenvalue $\lambda$ satisfies the estimate
  \begin{equation}
\int_{{\Bbb R}^{\rm d}} | \psi(x) |^2 e ^{2\delta |x|} dx < \ii
\label{eq:I1}\end{equation}
where $\delta$ depends on the distance $d(\lambda)=\dist \{\lambda,\sigma_{ess}(H)\}$
between this eigenvalue and the essential spectrum $\sigma_{ess}(H)$ of the operator $H$. More precisely, it was shown in \cite{Sh} that, for the Schr\"odinger operator $H=-\Delta+ V(x)$ with a semibounded real potential in the space $L_{2}({\Bbb R}^d)$, 
 estimate \e{eq:I1} holds with any $\delta$ such that
\[
\delta < 2^{-1}\ln(1+ d^2(\lambda)a^{-2})
\]
where $a$ depends only on 
\[
v_{0}=\inf_{x\in {\Bbb R}^d }V(x).
\]

Another method originated in papers  \cite{Combes, Connor} where eigenfunctions corresponding to
eigenvalues lying below $\sigma_{ess}(H)$ were considered for the $N$-particle 
Schr\"odinger operator $H$. In this case
estimate \e{eq:I1}   with any  
$\delta <\sqrt{d(\lambda)}$ was proven. 
 Then   Agmon \cite{Ag}   extended this result to very general differential operators of second order.
This estimate is of course more efficient and stronger than that of Shnol' but, in general,  does not hold for eigenvalues lying in gaps of $\sigma_{ess}(H)$. This can be shown on an example of an one-dimensional Schr\"odinger operator whose potential is a sum of a periodic function and of a  function with compact support.
 From a technical point of view,  an estimate of eigenfunctions for eigenvalues in gaps is essentially
more difficult than for eigenvalues below $\sigma_{ess}(H)$ because in the first case one has to work with the expression $\| Hf-\lambda f\|$ whereas it suffices to consider the quadratic form $(Hf,f)$  in the second case.

We note also paper  \cite{BCH} where an estimate of eigenfunctions for eigenvalues in gaps was obtained for   random Schr\"odinger operators.

Our goal here is to study exponential decay of eigenfunctions for  first order matrix differential operators
\begin{equation}
H =-i\sum_{j=1}^{\rm d} A_j \frac{\partial}{\partial x_{j}} +V(x)
\label{eq:I2}\end{equation}
acting in the space ${\cal H}=L_2({\Bbb R}^{\rm d} ; {\Bbb C}^n)$. Here 
$A_j=A_j^*$, $j=1,\ldots,  {\rm d}$, are constant matrices and $V(x)$ is a symmetric matrix-valued functiion.
 Set 
  \begin{equation}
\gamma =\max_{|x|=1} {\pmb|} \sum_{j=1}^{\rm d}  A_j   x_{j}{\pmb|}, \q 
{\pmb|} \cdot{\pmb|}={\pmb|} \cdot{\pmb|}_{{\Bbb C}^n}.
\label{eq:I3}\end{equation}
 For example, $\gamma=1$ for the Dirac operator. 
 Our main result is the estimate 
   \begin{equation}
\int_{{\Bbb R}^{\rm d}} {\pmb|} \psi(x) {\pmb|}^2 e ^{2\delta |x|} dx < \ii, \q
\forall   \: \delta< \gamma^{-1} d(\lambda),
\label{eq:A2}\end{equation}
  for all eigenvalues (including those lying in gaps of $\sigma_{ess}(H)$).  
   
 Of course for  first order   differential operators the problem considered is essentially simpler   than for operators of second order. This allows us to follow sufficiently closely the method of \cite{Ag} and to obtain estimate \e{eq:A2} without any specific assumptions on the coefficients $A_j$   and $V(x)$. Note that for second order elliptic operators one can pass from weighted $L_{2}$-estimate \e{eq:I1} to the pointwise estimate
% \begin{equation}
\[
|\psi(x,\lambda)|\leq C e^{-\delta |x|},\q \delta <\sqrt{d(\lambda)}.
\]
%\label{eq:I1a}\end{equation}
Such a  passage is eventually possible also for  first order systems but would require more or less stringent assumptions on coefficients of the system.

Here and below $C$ denotes   constants whose precise values are irrelevant; $I$ is the identity operator (in different spaces).

%***********************************************************
\section{The main result}
%*

% Let ${\Bbb B}_{R}$ be the ball $|x|\leq R$. 

% $f\in C_{0}^1 ({\Bbb R}^{\rm d}\setminus {\Bbb B}_{R}; {\Bbb C}^n)$

Let us first discuss the precise definition of the operator $H$.   Derivatives in \e{eq:I2} are  
understood in the sense of distributions.  We suppose that $V(x)$ is a locally bounded function, but we do not make any specific assumptions about its behavior at infinity. Instead of this we a priori assume that $H$ is selfadjoint on a domain with the following natural properties.

 \begin{assumption}\label{D}
Operator \e{eq:I2} is selfadjoint on   domain
${\cal D}(H)$ such that

$1^0$ If a (vector-valued) function $f (x)$ has compact support and
\[
 \sum_{j=1}^{\rm d} A_j  \partial f /\partial x_{j}\in L_{2}({\Bbb R}^{\rm d} ; {\Bbb C}^n),
 \]
 then $f \in {\cal D}(H)$.
 
 $2^0$ If $f\in {\cal D}(H)$, then
\[
 \sum_{j=1}^{\rm d} A_j  \partial f/\partial x_{j}\in L_{2}^{(loc)} ({\Bbb R}^{\rm d} ; {\Bbb C}^n).
 \]
\end{assumption}

If the function $V(x)$ is bounded, then the operator $H$ is selfadjoint on the  set of   functions   $f \in {\cal H}$ such that
\[
\int_{{\Bbb R}^{\rm d}} {\pmb|} \sum_{j=1}^{\rm d} A_j  \partial f/\partial x_{j}{\pmb|}^2 dx<\ii.
\]
In this case Assumption~\ref{D} is of course satisfied. Recall that the operator $H$ is called elliptic if the matrix 
  \begin{equation}
A(\xi)=     \sum_{j=1}^{\rm d} A_{j} \xi_{j}
\label{eq:A4}\end{equation}
 does not have zero eigenvalues for $\xi\neq 0$ or, equivalently,
\[
  {\pmb|} A(\xi) p {\pmb|}\geq c \, {\pmb|}  p {\pmb|} \, |\xi| , \q c>0, \q p \in {\Bbb C}^n,
   \q \xi\in {\Bbb R}^{\rm d}.
\]
Elliptic differential operators with bounded $V(x)$ are selfadjoint on  the Sobolev class $ \mathsf{H}^1 ({\Bbb R}^{\rm d} ; {\Bbb C}^n)=:{\cal D}(H)$.
 We emphasize however that we do not assume neither ellipticity of the operator $H$
 nor boundedness of $V(x)$.

 The following standard result is a  direct consequence of the spectral theorem.

 \begin{lemma}\label{1}
 For every $\varepsilon>0$ there exists $R$ such that 
    \begin{equation}
\| H f -\lambda f\| \geq (d(\lambda)-\varepsilon) \| f \|, \q
d(\lambda)=\dist \{\lambda,\sigma_{ess}(H)\},
\label{eq:AI}\end{equation}
for all $f\in   {\cal D}(H)$ such that $f(x)=0$ in the ball $|x|\leq R$.
       \end{lemma}
       
         \begin{pf}
         Suppose to the contrary. Then for some $\varepsilon>0$ and all $R$
          there exists $ f_{R}\in   {\cal D}(H)$, $\| f_{R} \|=1$, such that $f_{R}(x)=0$ for $|x|\leq R$ and
          \begin{equation}
\| H f_{R} -\lambda f_{R}\| < d(\lambda)-\varepsilon.
\label{eq:AIa}\end{equation}
Let us choose some $a< d(\lambda)$ and set  
$X=[\lambda- a, \lambda+ a]$,
$g_{R}=f_{R}- E(X) f_{R}$ where $E(X)$ is the spectral projection of the operator $H$ corresponding to the interval $X$. Since $g_{R}= E({\Bbb R}\setminus X) g_{R}$, the spectral theorem implies that
 \begin{equation}
\| H g_{R} -\lambda g_{R}\| >  a   \|  g_{R}\|.
\label{eq:AIb}\end{equation}
Using that $f_{R}\ri 0$ weakly and that the operator $E(X)$  has finite rank, we find that
\[
 \|  g_{R} -  f_{R}\| =o(1)  , \q \| H g_{R} - H f_{R}\| =o(1)
 \]
 as $R\ri\ii$. 
 Thus, comparing estimates \e{eq:AIa} and       \e{eq:AIb}, we arrive at the contradiction if 
 $ a> d(\lambda)- \varepsilon$ and $R$ is sufficiently large.
           \end{pf}

The following theorem is our main result.
 
  \begin{theorem}\label{2}
  Let the operator $H$ be defined by formula \e{eq:I2}, let Assumption~\ref{D} be satsified,  and let $\gamma$ be constant \e{eq:I3}.     If $H\psi=\lambda\psi$ where $\psi\in {\cal D}(H)$, then estimate   \e{eq:A2} holds.
       \end{theorem}
       
       \begin{pf}
   Put $ \varepsilon =2^{-1} (  d(\lambda) -  \gamma\delta)$ and  choose $R$   in such a way that  estimate  \e{eq:AI} holds if $f\in   {\cal D}(H)$ and  $f(x)=0$ for $|x|\leq R$.  Fix some number $R_{1}>R$.
       Let $\zeta\in C^\ii ({\Bbb R}_{+})$, $\zeta(r)=0$ for $r\leq R$ and $\zeta(r)=1$ for $r\geq R_{1}$.  For an arbitrary (large) $R_{2}> R_{1}$,  we choose $\rho \in C^\ii ({\Bbb R}_{+})$ such that $\rho(r)= \delta r$ for $r< R_{2} $, $\rho^\prime (r) \leq \delta  $ for  $r\in (R_{2},2R_{2})$ and $\rho(r)=  const$ for $r\geq 2 R_{2} $. Set $  \eta (r)= \zeta (r) e^{\rho(r)}$ 
   and $f (x) =\eta (|x|)\psi (x)$. Of course, $f =\psi-(1-\eta)\psi \in {\cal D}(H)$
   because $1-\eta \in C_{0}^\ii$. 
    It follows from equation $H\psi=\lambda\psi$ that
     %  \begin{equation}
     \[
(Hf)(x)- \lambda f(x)= -i r^{-1}    \eta^\prime (r) A(x) \psi (x)
\]
%\label{eq:A3}\end{equation} 
where the matrix $A(x)$ is defined by formula \e{eq:A4}.
Using notation \e{eq:I3}, we find that
 \begin{equation}
\| H f -\lambda f\|\leq \gamma \| \eta^\prime     \psi\|
 \leq \gamma(\| \rho^\prime e^\rho \zeta \psi\| + \|    e^\rho \zeta^\prime \psi \|).
\label{eq:A5}\end{equation}
 Since $\rho^\prime (r) \leq \delta  $, we have that
$$
\| \rho^\prime e^\rho \zeta \psi\| \leq  \delta \| f \|.
$$
Since $\supp \zeta^\prime \subset [R,R_{1}]$, we have that
$$
 \|    e^\rho \zeta^\prime \psi \| \leq C  \|     \psi \|
 $$
 where the constant $C$ does not depend on $R_{2}$.
Thus, it follows from  \e{eq:A5} that
\begin{equation}
\| H f -\lambda f\| \leq \gamma\delta \| f \| +  C  \|     \psi \|.
\label{eq:A5t}\end{equation}

On the other hand, the function $f$ satisfies estimate  \e{eq:AI} because $f(x)=0$ for $|x|\leq R$.
Comparing   estimates \e{eq:AI} where $ \varepsilon =2^{-1} (  d(\lambda) -  \gamma\delta)$  and \e{eq:A5t}, we see that
\[
(d(\lambda)-   \gamma\delta) \| f \| \leq C  \|     \psi \|
\]
and hence
\[
\int_{R_{1} \leq |x| \leq R_{2}} {\pmb|} \psi(x) {\pmb|}^2 e ^{2\delta |x|} dx \leq C  \|     \psi \|^2.
\]
This proves \e{eq:A2} because  $C$ does not depend on $R_{2}$
which can be chosen arbitrary large.
 \end{pf}
   
 %***********************************************************
\section{The Dirac operator}
%***********************************************************
 
 As an example to which Theorem~\ref{2}
directly applies, we
  now consider the Dirac operator describing a relativistic particle of
spin
$1/2$. Let  ${\cal H}=L_2({\Bbb R}^3; {\Bbb C}^4)$ and
\begin{equation}
 H = -i \sum_{j=1}^3 {\pmb \alpha}_j \frac{\partial}{\partial x_{j}} +m {\pmb \alpha}_0+
  W(x),
\label{eq:D1}\end{equation} 
where $m>0$ is the mass of  a particle and $4\times 4$ - Dirac
matrices  ${\pmb \alpha}_j={\pmb \alpha}_j^*$  
satisfy the anticommutation relations
$$
  {\pmb \alpha}_j {\pmb \alpha}_k +{\pmb \alpha}_k {\pmb \alpha}_j=0,\quad j\neq k,\quad 
{\pmb \alpha}_j^2=I,\quad  j,k=0,1,2,3.
$$ 
 These relations determine the matrices ${\pmb \alpha}_j$ up to a unitary equivalence in the space
${\Bbb C}^4$. For example, we can express them in terms of the Pauli matrices 
\begin{equation}
  \sigma_1=\left(\begin{array}{cc} 0& 1\\ 1& 0
\end{array}\right),
 \quad 
\sigma_2=\left(\begin{array}{cc} 0& -i\\
 i &0
\end{array}\right),  \quad 
\sigma_3=\left(\begin{array}{cc} 1& 0\\
 0 & -1
\end{array}\right)
\label{eq:D2}\end{equation}
by the relation
 \begin{equation}
  {\pmb \alpha}_i=\left(\begin{array}{cc} 0& \sigma_i\\
\sigma_i& 0
\end{array}\right),
\quad i=1,2,3,\quad 
{\pmb \alpha}_0=\left(\begin{array}{cc} I& 0\\
 0& -I
\end{array}\right).
\label{eq:D3}\end{equation} 
We do not require any special assumptions on the matrix-valued function $W(x)$.

Operator \e{eq:D1} has form \e{eq:I2} where $d=3$, $n=4$, $A_{j}={\pmb \alpha}_j$ and
$$
V(x)=m {\pmb \alpha}_0+  W(x). 
$$
It follows from formulas \e{eq:D2}, \e{eq:D3} that  matrix \e{eq:A4} now equals
% \begin{equation}
\[
  A(\xi)=\left(\begin{array}{cc} 0& \sum_{j=1}^3\sigma_j \xi_{j}\\
\sum_{j=1}^3\sigma_j \xi_{j} & 0
\end{array}\right).
\]
% \label{eq:D4}\end{equation} 
Since $(\sum_{j=1}^3\sigma_j \xi_{j})^2=|\xi|^2$,
an easy calculation shows that the matrix $A(\xi)$ has two eigenvalues $|\xi | $ and $- |\xi | $ (of multiplicity $2$ each) and hence expression \e{eq:I3}  equals $1$. Thus, according to Theorem~\ref{2} estimate \e{eq:A2} where $\gamma=1$ is true for all eigenfunctions of the Dirac operator. Perhaps this result is new even for perturbations $W(x)$ decaying at infinity when 
$\sigma_{ess}(H)=(-\ii, -m^2]\cup [m^2, \ii)$ has the gap $ (  -m^2 , m^2)$.

Recall that estimate  \e{eq:I1} with $\delta <\sqrt{d(\lambda)}$ is sharp for the Schr\"odinger operator. We cannot claim the same about our estimate for the Dirac operator.
To clarify these statements, let us consider the simplest one-dimensional case. 

Let   first $H=-d^2/dx^2 + V(x)$  where $V(x)$ is a function with compact support. If $\lambda<0$ is an eigenvalue of $H$, then, for large $\pm x >0 $, the corresponding eigenfunction equals $\psi(x)=c_{\pm}e^{-\sqrt{|\lambda|} |x|}$   where   $c_{\pm}$ are some constants. Since now $d(\lambda)=|\lambda|$, one cannot hope to prove estimate
 \e{eq:I1} even for $\delta  =\sqrt{d(\lambda)}$.
 
 The one-dimensional   Dirac operator is defined by the formula
\[
 H = -i \left(\begin{array}{cc} 0& 1\\ 1& 0
\end{array}\right) d/dx+ \left(\begin{array}{cc} m& 0\\ 0& -m
\end{array}\right)+  W(x)
\]
in the space ${\cal H}=L_2({\Bbb R}; {\Bbb C}^2)$. If the (matrix-valued) function $W(x)$
has compact support and  $\lambda\in (  -m^2 , m^2)$ is an eigenvalue of the operator $H$, then,   for large $|x|$,  the corresponding eigenfunction equals 
$$
\psi(x)=c_{\pm} (\sqrt{m+\lambda},
i\sqrt{m-\lambda} \sgn x )^t   e^{-\sqrt{m^2-\lambda^2} |x|}, \q \pm x>0.
$$
    Thus, the integral in  \e{eq:A2} is convergent for all $\delta <\sqrt{m^2-\lambda^2}$
   which is better than $\delta  < d(\lambda)=m-|\lambda|$ as required  in  \e{eq:A2}.

 %***********************************************************
\section{The Hill operator}
%***

First we discuss a curious application of the Agmon estimate to the one-dimensional Schr\"odinger operator 
   \begin{equation}
 H=-d^2/dx^2 + V(x)
\label{eq:Hill}\end{equation}
   with a periodic potential $V(x)=V(x+1)$. Recall that for any $\lambda$ in the gap of the spectrum $\sigma(H)=\sigma_{ess}(H)$ of $H$ the equation
\begin{equation}
-y^{\prime \prime}(x)+ V(x)y(x)=\lambda y(x)
\label{eq:H0}\end{equation}
has two non-trivial solutions $y_{\pm}(x)$ admitting the representation
  \begin{equation}
y_{\pm}(x)=\rho^{\mp x} p_{\pm}(x)
\label{eq:H}\end{equation}
where both functions $p_{\pm}(x) $ are either periodic or antiperiodic 
($p_{\pm}(x)=  p_{\pm}(x+1)$ or $p_{\pm}(x)= - p_{\pm}(x+1)$) and the number $\rho=\rho(\lambda)> 1$ is known as the multiplicator. Thus, the solution $y_{\pm}(x)$ exponentially decays as $x\ri \pm \ii$.

\begin{proposition}\label{H}
Suppose that
\[
\lambda< \lambda_{0}:=\inf \sigma_{ess}(H).
\]
Then
   \begin{equation}
 \ln^2 \rho(\lambda)\geq \lambda_{0}-\lambda.
\label{eq:H1}\end{equation}
       \end{proposition}
       
       \begin{pf}
       Let us consider the auxiliary family of operators
          \begin{equation}
          H_{\alpha}=-d^2/dx^2 + V(x)-\alpha Q(x)
          \label{eq:H2}\end{equation}
        where $Q \in C_{0}^\ii ({\Bbb R})$ and $Q (x)\geq 0$. If $Q (x)$ is not identically zero, then we can choose $\alpha > 0$ such that $\lambda$ is an eigenvalue of the operator $H_{\alpha}$. The corresponding eigenfunction
         \begin{equation}
          \psi(x)=c_{\pm}y_{\pm}(x)
          \label{eq:H2a}\end{equation}
           for sufficiently large $\pm x > 0$. On the other hand, it satisfies estimate \e{eq:I1} where $\delta$ is any number smaller than        $\sqrt{\lambda_{0}-\lambda}$. By virtue of \e{eq:H} this is possible only if estimate          \e{eq:H1} is true.
        \end{pf}
        
         Finally, we show that estimate \e{eq:I1}
   with   arbitrary $\delta< \sqrt{d(\lambda)}$ is not in general true for eigenvalues in gaps. To that end,  we need some additional, although always quite elementary, information on  the Hill operator. Let $\varphi(x,\lambda)$ and $\theta (x,\lambda)$ be the solutions of equation \e{eq:H0} satisfying the conditions
   $$
   \varphi(0,\lambda)=0,\q \varphi^\prime(0,\lambda)=1,\q  \theta(0,\lambda)=1,\q 
   \theta^\prime(0,\lambda)=0,
   $$ 
   and let
   $$
   F(\lambda)=2^{-1}(\varphi^\prime(1,\lambda)+ \theta(1,\lambda)).
   $$
   Then $|F(\lambda)| > 1$ for regular points of $H$, and the equation $| F(\lambda) |=  1$ determines   band edges of its spectrum. The functions $\rho(\lambda)$  and $F(\lambda)$ are related by the formula
\begin{equation}
  \rho(\lambda) =| F(\lambda) | +\sqrt{F^2(\lambda) -1} .
  \label{eq:H3}\end{equation}
We note that the functions $p_{\pm}(x)$ in   \e{eq:H} are periodic (antiperiodic) if 
$F(\lambda)  > 1$ (if  $F(\lambda)  < - 1$).
  
  Let $\lambda$ be a regular point close to one of two (one for $\lambda< \sigma (H)$)
  band edges. Denote this edge by  $\lambda_{*}$.   It follows from \e{eq:H3} that,  as $\lambda\ri \lambda_{*}$, 
  \[
  \rho(\lambda)=1 + \sqrt{2 |F^\prime(\lambda_{*}) (\lambda-\lambda_{*})|}+ O(\lambda-\lambda_{*})
  \]
   and therefore  
   \begin{equation}
    \ln^2 \rho(\lambda)  = 2 |F^\prime(\lambda_{*})  |\, 
    | \lambda-\lambda_{*}|+ O( |\lambda-\lambda_{*}|^{3/2}).
     \label{eq:H4}\end{equation}

     Let us now use standard asymptotic formulas
     \[
     \varphi^\prime(1,\lambda)= \cos \sqrt{\lambda}+ O\Big(\frac{1}{\sqrt{\lambda}}\Big),\q 
      \theta(1,\lambda)=\cos \sqrt{\lambda}+ O\Big(\frac{1}{\sqrt{\lambda}}\Big), \q \lambda\ri \ii,
     \]
     which can  easily be deduced from the Volterra integral equations for solutions
     $ \varphi (x,\lambda)$ and $ \theta (x,\lambda)$     of differential equation \e{eq:H0}. Moreover, these asymptotic formula can be differentiated so that
      \begin{equation}
F^\prime(\lambda)= -\frac{\sin \sqrt{\lambda}}{2\sqrt{\lambda}}
+ O\Big(\frac{1}{ \lambda}\Big), \q \lambda\ri \ii.
  \label{eq:H5}\end{equation}
  
  Comparing formulas \e{eq:H4} and \e{eq:H5}, we obtain the following result.
  
  \begin{proposition}\label{H1}
  For an arbitrary $c>1$, sufficiently distant gaps and $\lambda$ sufficiently close to a band edge $ \lambda_{*}$   the estimate
% \begin{equation}
\[
    \ln^2 \rho(\lambda)  \leq c \lambda_{*}^{-1/2}
    | \lambda-\lambda_{*}|  
    \]
%         \label{eq:H6}\end{equation}
         is true.
       \end{proposition}

     Let us again consider operator \e{eq:H2}.  We use the following simple result.
     
       \begin{lemma}\label{H2}
     For an arbitrary regular point     $\lambda\in {\Bbb R}$ of the operator $H$, there exists $Q =G^2 \in C_{0}^\ii ({\Bbb R})$ and  $\alpha\in {\Bbb R}$ such that the operator  $H_{\alpha}$ has the eigenvalue $\lambda$.
     \end{lemma}
     
    \begin{pf} 
     It suffices to show that, for suitable $G$ and $\alpha$,  the Birman-Schwinger equation
      \begin{equation}
   u-\alpha G (H-\lambda I)^{-1} G u=0, \q H=H_{0} ,
     \label{eq:H7}\end{equation}
     has a non-trivial solution $u$. Let us check first that we can find $G$ such that the operator 
     $G (H-\lambda I)^{-1} G\neq 0$. Supposing to the contrary, we denote by $G_{n}$ the operator of multiplication by the function $G_{n}(x)=g(x/n)$ where $g \in C_{0}^\ii ({\Bbb R})$ and $g(x)=1$ for $|x|\leq 1$. By our assumption,   $G_{n} (H-\lambda I)^{-1} G_{n}= 0$ and hence
     \[
     (G_{n} (H-\lambda I)^{-1} G_{n} u_{1}, u_{2})=0
     \]
     for all $n$ and all $u_{1}, u_{2}\in L_{2}({\Bbb R})$. Since $G_{n}$ converges strongly to the identity operator, we have that
     $   ((H-\lambda I)^{-1}   u_{1}, u_{2})=0$. Since $u_{2}$ is arbitrary, this implies that
    $ (H-\lambda I)^{-1}   u_{1}=0$ and hence $u_{1}=0$ although $u_{1} $ is arbitrary.
    
    Let us pick   an arbitrary eigenvalue $\mu\neq 0$ of the compact operator $G (H-\lambda I)^{-1} G$. Then for $\alpha =\mu^{-1}$, the corresponding eigenfunction $u$ satisfies equation \e{eq:H7}.
        \end{pf} 
  
   Now we can construct a counterexample to estimate \e{eq:I1}
   with   arbitrary $\delta< \sqrt{d(\lambda)}$  for eigenvalues in gaps.   
   
     \begin{counterexample}\label{H5}
Let a periodic function $V(x)$ be such that        operator  \e{eq:Hill}  has  infinite number of gaps.       Then for an arbitrary $\varepsilon>0$, there exists an eigenvalue $\lambda$ of operator \e{eq:H2}      where $Q   \in C_{0}^\ii ({\Bbb R})$ such that the corresponding eigenfunction equals
       \begin{equation}
    \psi(x)=  e^{- \delta |x|} p_{\pm}(x)
    \label{eq:H8}\end{equation}
    with $\delta < \varepsilon  \sqrt{d(\lambda)}$ for sufficiently large $\pm x >0$.
    Here both $p_{\pm}(x)$ are either periodic or antiperiodic functions which are not identically zero.
  \end{counterexample}
  
  Indeed, using Proposition~\ref{H1}, we can find $\lambda$ such that $\ln \rho(\lambda)<\varepsilon  \sqrt{d(\lambda)}$. Then, by Lemma~\ref{H2}, we can construct $Q(x)$ such that operator \e{eq:H2} has eigenvalue $\lambda$. The corresponding eigenfunction 
  $\psi(x)$ is given by formula  \e{eq:H2a}        for sufficiently large $\pm x > 0$. By virtue of   \e{eq:H} this leads to representation \e{eq:H8} with $\delta = \ln \rho(\lambda)$.

  It is possible that the results of this section can be deduced from very general results of \cite{Ag1}.
    
        Thanks are due to N. Firsova for a useful discussion of the periodic problem.

      \end{document}